\documentclass [reqno] {conm-p-l} 

\def\newpic#1{%
   \def\emline##1##2##3##4##5##6{%
      \put(##1,##2){\special{em:point #1##3}}
      \put(##4,##5){\special{em:point #1##6}}
      \special{em:line #1##3,#1##6}}}
\newpic{}

\copyrightinfo{2000}{American Mathematical Society}

\newtheorem{theorem}{Theorem}[section]
\newtheorem{lemma}[theorem]{Lemma}
\newtheorem{corollary}[theorem]{Corollary}

\theoremstyle{definition}
\newtheorem{definition}[theorem]{Definition}

\theoremstyle{remark}

\numberwithin{equation}{section}

\begin{document}

\title{Equations of Maxwell Type}

\author{K. O. Makhmudov}

\address{Department of Mechanics and Mathematics,
         University of Samarkand,
         University Boulevard 15,
         703004 Samarkand,
         Uzbekistan}

\email{komil.84@mail.ru}

\author{O. I. Makhmudov}

\address{Department of Mechanics and Mathematics,
         University of Samarkand,
         University Boulevard 15,
         703004 Samarkand,
         Uzbekistan}

\email{olimjan@yahoo.com}

\author{N. Tarkhanov}

\address{Institute of Mathematics,
         University of Potsdam,
         Am Neuen Palais 10,
         14469 Potsdam,
         Germany}

\email{tarkhanov@math.uni-potsdam.de}

\thanks{The authors wish to express their thanks to V.~Kravchenko for several
        helpful comments concerning boundary value problems for the classical
        Maxwell system.
        The second author gratefully acknowledges the financial support of
        the Deutscher Akademischer Austauschdienst.}

\date{September 2, 2009}


\subjclass [2000] {Primary 35Q60; Secondary 35Cxx}

\keywords{Electromagnetic waves,
          scattering,
          elliptic complex,
          Green formulas,
          Stratton-Chu formulas,
          Cauchy problem}

\begin{abstract}
For an elliptic complex of first order differential operators on a smooth
manifold $\mathcal{X}$, we define a system of two equations which can be
thought of as abstract Maxwell equations.
The formal theory of this system proves to be very similar to that of
classical Maxwell's equations.
The paper focuses on boundary value problems for the abstract Maxwell
equations, especially on the Cauchy problem.
\end{abstract}

\maketitle

\tableofcontents

\section*{Introduction}
\label{s.I}

The term Maxwell's equations applies to a set of eight equations published by
Maxwell in \cite{Maxw65}.
He called them ``general equations of the electromagnetic field.''

To be more specific, we will restrict ourselves to reduced Maxwell's equations
which describe the propagation of electromagnetic waves in a homogeneous
isotropic medium in $\mathbb{R}^3$ with
   electric permittivity $\varepsilon$,
   magnetic permeability $\mu$
and
   electrical conductivity $\sigma$.
An electromagnetic wave is described by the vectors of
   electromotive force $E (t,x)$
and
   magnetic force $H (t,x)$.
Then, Maxwell's equations are
\begin{equation}
\label{eq.Maxwell}
\begin{array}{rcl}
   - \varepsilon\, (d/dt) E + \mathrm{curl}\, H
 & =
 & \sigma E,
\\
   \mu\, (d/dt) H + \mathrm{curl}\, E
 & =
 & 0,
\end{array}
\end{equation}
see for instance \cite{Jack75}.

Maxwell's equations have a close relation to special relativity.
Not only were Maxwell's equations a crucial part of the historical development
of special relativity, but also special relativity has motivated a compact
mathematical formulation of Maxwell's equations.

Although Maxwell's equations apply throughout space and time, practical
problems are finite and solutions to Maxwell's equations inside the solution
region $\mathcal{X}$ are joined to the remainder of the universe through
boundary conditions, cf. \cite{Monk03},
                         \cite{SeniVola95},
                         \cite{Frie90},
                         etc.,
and started in time using initial conditions, cf. for instance
   \cite{HarmHuss94}.

In \cite{MakhNiyoTark08} we constructed an explicit formula which restores
solutions of the oscillation equation of the couple-stress theory of
elasticity in a bounded domain $\mathcal{X}$ in $\mathbb{R}^3$ through the
given displacement and stress values on a part $\mathcal{S}$ of the boundary
of the domain.
The main ingredient of our construction is an expansion of the fundamental
solution $\mathit{\Phi} (x-y)$ of couple-stress elasticity.
The equations of couple-stress elasticity factorise the Helmholtz equation
in $\mathbb{R}^3$.
It follows that $\mathit{\Phi} (x-y)$ amounts to the quotient applied to the
fundamental solution of the Helmholtz equation.
In order to expand this latter we used spherical harmonics in $\mathbb{R}^3$
and Bessel functions.

Maxwell's equations factorise the Helmholtz equation in $\mathbb{R}^3$, too,
and so the construction of \cite{MakhNiyoTark08} goes through in this case in
much the same way.
This work was intended as an attempt at constructing an explicit formula which
restores solutions of Maxwell's equations in a bounded domain $\mathcal{X}$ in
$\mathbb{R}^3$ through their values on a part $\mathcal{S}$ of the boundary of
$\mathcal{X}$.
While such a formula is of practical interest and, to our best knowledge, new,
it would not be surprising, for the arguments are routine.

However, on thinking of the so-called Stratton-Chu formulas
   \cite{Stra41}
we have revealed a very rich algebraic structure which goes beyond classical
Maxwell's equations.
Eliminating the electric or magnetic field from Maxwell's equations leads to
the Helmholtz equation for the components of the other field.
Algebraically the elimination means multiplication of Maxwell's operator with
a matrix of scalar differential operators from the left.
Any such matrix is referred to as the quotient of the Helmholtz and Maxwell's
operators.
In contrast to Dirac type operators which multiply by their formal adjoints to
the Laplace operator, the quotient of the Helmholtz and Maxwell's operators is
of order $2$.

Hence it follows that the potential theory for Maxwell's equations is not
quite standard.
To handle simple and double layer potentials related to these equations one
has to take care of appropriate modification of the results in
   \cite[Ch.~2]{Tark90}.
The nature of such modification becomes most transparent if one specifies
Maxwell's equations within a wider class of differential equations associated
with elliptic complexes.

For this purpose, we recall in Section \ref{s.aMe} the concept of Maxwell's
equations.
When written in terms of differential forms, they can readily be interpreted
in the context of arbitrary elliptic complexes.
The equations obtained in this way we call abstract Maxwell equations, and
they raise new boundary value problems.
In Section \ref{s.SCf} we generalise the Stratton-Chu formula to solutions of
abstract Maxwell equations.
This allows one to construct explicitly the so-called Calderon projections,
   see Theorem 5.4 of \cite[1.4.1]{PikeSaba02},
which reduce the study of boundary value problems to boundary integral
equations,
   cf. Section \ref{s.bvp}.
The Cauchy problem with data on the whole boundary for classical Maxwell's
equations is certainly overdetermined.
In \cite{Krav92} and
   \cite{KhmeKravRabi03},
some solvability criteria for this problem are proved by using biquaternionic
functions.
The Cauchy problem with data on a nonempty part $\mathcal{S}$ of the boundary
is of great importance for analysts.
We discuss this problem in Sections \ref{s.tCp},
                                    \ref{s.Eotfs} and
                                    \ref{s.reg}.
The method of \cite{YarmIshaMakh92} and
              \cite{Makh04}
applies to construct a Carleman function for the Helmholtz operator, and so
for all its matrix factorisations, in the case of conical domains.
Here, we use another approach which is elaborated in
   \cite{Tark95} and
   \cite{Shla96}.
In Section \ref{s.poemw} we apply these results to classical Maxwell's
equations.

\section{Abstract Maxwell's equations}
\label{s.aMe}

For the sake of simplicity we confine ourselves to the case of time-harmonic
electromagnetic waves.
Then the electric and magnetic fields are of the asymptotic form
$$
\begin{array}{rcl}
   E (t,x)
 & =
 & (\varepsilon + \imath\, \sigma / \omega)^{- 1/2}
   e^{- \imath \omega t} E (x),
\\
   H (t,x)
 & =
 & \mu^{- 1/2}
   e^{- \imath \omega t} H (x).
\end{array}
$$
From Maxwell's equations (\ref{eq.Maxwell}) with arbitrary time dependence we
deduce easily that the parts of $E$ and
                                $H$
depending on the space coordinates satisfy Maxwell's equations
$$
\begin{array}{rcl}
   \imath k\, E + \mathrm{curl}\, H
 & =
 & 0,
\\
   - \imath k\, H + \mathrm{curl}\, E
 & =
 & 0,
\end{array}
$$
where $k$ is the wave constant defined by
$
   k^2 = (\varepsilon + \imath\, \sigma / \omega) \mu \omega^2.
$
The sign of $k$ is chosen from the condition $\Im k \geq 0$.

The expression of Maxwell's equations in terms of differential forms leads to
a further notational and conceptual simplification.
On introducing the de Rham complex in $\mathbb{R}^3$
$$
   0
 \longrightarrow
   \mathit{\Omega}^0 (\mathbb{R}^3)
 \stackrel{d}{\longrightarrow}
   \mathit{\Omega}^1 (\mathbb{R}^3)
 \stackrel{d}{\longrightarrow}
   \mathit{\Omega}^2 (\mathbb{R}^3)
 \stackrel{d}{\longrightarrow}
   \mathit{\Omega}^3 (\mathbb{R}^3)
 \longrightarrow
   0
$$
we can think of $E$ as a differential form $u$ of degree $1$,
                $H$ as a differential form $f$ of degree $2$,
thus identifying $\mathrm{curl}\, E$ with $du$ and
                 $\mathrm{curl}\, H$ with $d^\ast f$.
Here, $d^\ast$ stands for the formal adjoint operator of $d$.
In this way Maxwell's equations can be written in the form
\begin{equation}
\label{eq.space}
\begin{array}{rcl}
   \imath k\, u + d^\ast f
 & =
 & 0,
\\
   - \imath k\, f + du
 & =
 & 0,
\end{array}
\end{equation}
which already make sense not only for differential forms $u$ and
                                                         $f$
of degree $1$ and
          $2$
in $\mathbb{R}^3$, respectively, but also for differential forms $u$ and
                                                                 $f$
of degree $i$ and
          $i+1$
in $\mathbb{R}^n$, where $-1 \leq i \leq n$.

For a recent account of mathematical treatment of the electromagnetic
scattering theory we refer the reader to
   \cite{ColtKres98},
   \cite{PikeSaba02},
etc.

Equations (\ref{eq.space}) generalise to arbitrary complexes of differential
operators on a $C^\infty$ manifold $\mathcal{X}$ with or without boundary.
More precisely, consider a complex of first order differential operators on
$\mathcal{X}$ acting in sections of vector bundles over $\mathcal{X}$, i.e.
\begin{equation}
\label{eq.complex}
   0
 \longrightarrow
   C^{\infty} (\mathcal{X},F^0)
 \stackrel{A^0}{\longrightarrow}
   C^{\infty} (\mathcal{X},F^1)
 \stackrel{A^1}{\longrightarrow}
   \ldots
 \stackrel{A^{N-1}}{\longrightarrow}
   C^{\infty} (\mathcal{X},F^N)
 \longrightarrow
   0
\end{equation}
where
   $A^i \in \mathrm{Diff}^1 (F^i, F^{i+1})$ satisfy $A^{i+1} A^i = 0$ for all
   $i$.
As usual,
   we write $A^i u$ simply $Au$ for $u \in C^{\infty} (\mathcal{X},F^i)$,
when no confusion can arise.

By $F^i$ is meant a smooth vector bundle of rank $k_i$ which is nonzero only
for $i = 0, 1, \ldots, N$.
We give $F^i$ a Hermitean metric, i.e. scalar products
   $x \mapsto (v,w)_x$
in the fibres $F^i_x$ which smoothly depend on the point $x \in \mathcal{X}$.
This defines a conjugate linear isomorphism $\ast$ of $F^i$
                         to the algebraic dual bundle $F^i{}^\ast$
by
   $\langle \ast w, v \rangle_x = (v,w)_x$
for $v, w \in F^i_x$.

Fix a smooth positive volume form $dx$ on $\mathcal{X}$.
This yields a scalar product on $C^{\infty}_{\mathrm{comp}} (\mathcal{X},F^i)$
by
$$
   (u,v)
 = \int_{\mathcal{X}} (u (x),v(x))_x dx
$$
for $u, v \in C^{\infty}_{\mathrm{comp}} (\mathcal{X},F^i)$.
The completion of this space with respect to the corresponding norm is denoted
by $L^2 (\mathcal{X},F^i)$.
We moreover introduce the formal adjoint
   $A^i{}^\ast \in \mathrm{Diff}^1 (F^{i+1},F^i)$
for each operator $A^i$ by requiring
   $(Au,g) = (u,A^\ast g)$
for all $u \in C^{\infty} (\mathcal{X},F^i)$ and
        $g \in C^{\infty} (\mathcal{X},F^{i+1})$
whose supports do not meet each other on the boundary of $\mathcal{X}$.

The formal selfadjoint operators
   $\mathit{\Delta}^i = A^i{}^\ast A^i + A^{i-1} A^{i-1}{}^\ast$
on $\mathcal{X}$ are called the Laplacians of (\ref{eq.complex}).
The ellipticity of this complex joust amounts to the fact that each Laplacian
$\mathit{\Delta}^i \in \mathrm{Diff}^2 (F^i)$ is an elliptic operator of order
two,
   cf. \cite[2.1.4]{Tark90}.

Given any $-1 \leq i \leq N$,
   by Maxwell's equations for complex (\ref{eq.complex}) at step $i$
are meant
\begin{equation}
\label{eq.abstract}
\begin{array}{rcl}
   - \varepsilon\, (d/dt) u + A^\ast f
 & =
 & \sigma u,
\\
   \mu\, (d/dt) f + Au
 & =
 & 0,
\end{array}
\end{equation}
with $u$ and
     $f$
being unknown functions of $t$ with values in sections of $F^i$ and
                                                          $F^{i+1}$,
respectively.
When looking for time-harmonic solutions to these equations of the form
$$
\begin{array}{rcl}
   u (t,x)
 & =
 & (\varepsilon + \imath\, \sigma / \omega)^{- 1/2}
   e^{- \imath \omega t} u (x),
\\
   f (t,x)
 & =
 & \mu^{- 1/2}
   e^{- \imath \omega t} f (x),
\end{array}
$$
one arrives at stationary equations for the parts of $u$ and $f$ that depend
only on the points of the underlying manifold $\mathcal{X}$.
These are
$$
\begin{array}{rcl}
   \imath k\, u + A^\ast f
 & =
 & 0,
\\
   - \imath k\, f + Au
 & =
 & 0,
\end{array}
$$
where $k$ is the wave constant defined above.

\begin{definition}
\label{d.Maxwell}
Let $-1 \leq i \leq N$.
By the Maxwell operator for complex (\ref{eq.complex}) at step $i$ is meant
$$
   M^i
 = \Big( \begin{array}{cc}
           \imath k
         & A^i{}^\ast
\\
           A^i
         & - \imath k
         \end{array}
   \Big).
$$
\end{definition}

As is usual in homological algebra, we will omit the index $i$ of $M^i$ when
it is clear from the context.

By definition, $M^i$ is a first order differential operator from sections of
$F^i \oplus F^{i+1}$ to sections of the same bundle over $\mathcal{X}$.
This operator fails to be elliptic of order $1$ in the classical sense unless
$N = 2$.
On the other hand,
   applying $A^\ast$ to both sides of $\imath k\, u + A^\ast f = 0$
we conclude that $A^{i-1}{}^\ast u = 0$ unless $k = 0$.
Analogously, from $- \imath k\, f + Au = 0$ it follows that $A^{i+1} f = 0$
unless $k = 0$.
Complementing Maxwell's equations by their differential consequences
   $A^{i-1}{}^\ast u = 0$ and
   $A^{i+1} f = 0$
yields a system of first order differential equations for $u$ and
                                                          $f$,
whose classical symbol is injective.
Another way of stating this is to say that there is a differential operator
$C^i$ from sections of $F^i \oplus F^{i+1}$ to sections of the same bundle,
such that $C^i M^i$ is a second order differential operator on $\mathcal{X}$
elliptic in the classical case.
An easy computation shows that
\begin{equation}
\label{eq.quotient}
   C^i
 = \Big(
         \begin{array}{cc}
           \imath k + (1/\imath k) A^{i-1} A^{i-1}{}^\ast
         & A^i{}^\ast
\\
           A^i
         & - \imath k - (1/\imath k) A^{i+1}{}^\ast A^{i+1}
         \end{array}
   \Big).
\end{equation}

\begin{lemma}
\label{l.factorisation}
As defined above, $C^i$ satisfies
$$
   C^i M^i
 = M^i C^i
 = \Big(
         \begin{array}{cc}
           \mathit{\Delta}^i - k^2
         & 0
\\
           0
         & \mathit{\Delta}^{i+1} - k^2
         \end{array}
   \Big).
$$
\end{lemma}

\begin{proof}
This is straightforward.
\end{proof}

\section{Stratton-Chu formula}
\label{s.SCf}

For classical Maxwell's equations, the so-called Stratton-Chu formula
   \cite{Stra41}
lies in the base of application of analytical methods, cf.
   \cite[4.2]{ColtKres98}.
From the point of view of contemporary analysis this is just a very particular
case of Green formulas, cf.
   \cite[2.5.4]{Tark90}.

Let
   $\mathcal{X}$ be a compact manifold with boundary smoothly embedded into a
   larger $C^\infty$ manifold $\mathcal{X}'$,
and
   (\ref{eq.complex}) be defined on all of $\mathcal{X}'$.
Suppose that both $\mathit{\Delta}^i - k^2$ and
                  $\mathit{\Delta}^{i+1} - k^2$
are elliptic of order two and satisfy the uniqueness condition for the local
Cauchy problem on $\mathcal{X}'$.
Then these operators have left fundamental solutions on $\mathcal{X}'$ which
we denote by $G^i$ and
             $G^{i+1}$,
respectively.
These latter are classical pseudodifferential operators of order $-2$ acting
on sections of vector bundles $F^i$ and
                              $F^{i+1}$
over $\mathcal{X}'$.

\begin{lemma}
\label{l.fs}
The pseudodifferential operator
$$
   \mathit{\Phi}^i
 = \Big(
         \begin{array}{cc}
           G^i\, (\imath k + (1/\imath k) A^{i-1} A^{i-1}{}^\ast)
         & G^i\, A^i{}^\ast
\\
           G^{i+1}\, A^i
         & G^{i+1}\, (- \imath k - (1/\imath k) A^{i+1}{}^\ast A^{i+1})
         \end{array}
   \Big)
$$
is a left fundamental solution of the Maxwell operator $M^i$ on $\mathcal{X}'$.
\end{lemma}

\begin{proof}
From Lemma \ref{l.factorisation} it follows immediately that
$$
   \mathit{\Phi}^i
 = 
   \Big(
         \begin{array}{cc}
           G^i
         & 0
\\
           0
         & G^{i+1}
         \end{array}
   \Big)
 \circ
   C^i
$$
is a left fundamental solution of $M^i$.
It remains to substitute the explicit expression (\ref{eq.quotient}) for $C^i$.
\end{proof}

Note that $\mathit{\Phi}^i$ is a pseudodifferential operator of order $0$ (not
$-1$) on sections of vector bundle $F^i \oplus F^{i+1}$ over $\mathcal{X}'$.
This corresponds to the fact that the Maxwell operator is elliptic in the sense
of Douglis-Nirenberg with weights $s_i = t_j = 0$, provided $k \neq 0$.
As usual, we use the same letter to designate
   the operator $\mathit{\Phi}^i$
acting
   from smooth sections of compact support of the bundle $F^i \oplus F^{i+1}$
   to distribution sections of the same bundle,
and
   its Schwartz kernel on $\mathcal{X}' \times \mathcal{X}'$ which is denoted
   by $\mathit{\Phi}^i (x,y)$.

\begin{theorem}
\label{t.Green}
For each
   $u \in H^1 (\mathcal{X}, F^i)$ and
   $f \in H^1 (\mathcal{X}, F^{i+1})$
satisfying $M^i (u,f) = 0$ in the interior of $\mathcal{X}$, it follows that
$$
     \int_{\partial \mathcal{X}}
     \mathit{\Phi}^i (x,y)
     \binom {\sigma (A^i)^\ast (y,\imath \nu (y)) f (y)}
            {\sigma (A^i) (y,\imath \nu (y)) u (y)}
     ds (y)
 =
   \left\{ \begin{array}{rcl}
             \displaystyle
             \binom {u (x)}{f (x)},
           & \mbox{if}
           & x \in \mathcal{X} \setminus \partial \mathcal{X},
\\
             0,
           & \mbox{if}
           & x \in \mathcal{X}' \setminus \mathcal{X}.
           \end{array}
   \right.
$$
\end{theorem}

Here,
   $\nu (y)$ is the outward unit normal vector to the boundary at a point
   $y \in \partial \mathcal{X}$,
   $ds$ the surface measure on $\mathcal{X}$,
and
   $\sigma (A^i)$ the (classical) symbol of $A^i$.

\begin{proof}
This formula is a particular case of the general Green formula
   \cite[2.5.4]{Tark90},
for
$$
   \sigma (M^i)
 = \Big(
         \begin{array}{cc}
           0
         & \sigma (A^i)^\ast
\\
           \sigma (A^i)
         & 0
         \end{array}
   \Big).
$$
\end{proof}

Write $t (u)$ and
      $n (f)$
for the tangential part of $u$ and
    the normal part of $f$
on the boundary of $\mathcal{X}$ relative to complex (\ref{eq.complex}),
   respectively,
see \cite[3.2.2]{Tark90}.
Then we obtain
\begin{equation}
\label{eq.decomposition}
     \binom {\sigma (A^i)^\ast (\cdot,\imath \nu) f}
            {\sigma (A^i) (\cdot,\imath \nu) u}
 =
     \binom {\imath \sigma (\mathit{\Delta}^i) (\cdot,\nu) n (f)}
            {\imath \sigma (A^i) (\cdot,\nu) t (u)}
\end{equation}
on $\partial \mathcal{X}$.
A familiar argument shows that
   Theorem \ref{t.Green} still holds for all solutions $(u,f)$ of Maxwell's
   equations in the interior of $\mathcal{X}$,
such that $t (u)$ and
          $n (f)$
have weak limit values on the boundary $\partial \mathcal{X}$,
   cf. \cite[9.4]{Tark95}.

For classical Maxwell's equations the formula of Theorem \ref{t.Green} is
known as the Stratton-Chu formula, see \cite{Stra41}.
In that case it manifests a more refined structure,
   for the fundamental solution $\mathit{\Phi}^i$ can be written explicitly.
This can indeed be done in the context of the so-called complexes of Dirac
type in $\mathbb{R}^n$ which are characterised by the property that all the
Laplacians
   $\mathit{\Delta}^i$
are diagonal operators, with the usual (nonnegative) scalar Laplace operator
on the diagonal.
In other words,
$
   \mathit{\Delta}^i
 = - E_{k_i} \mathit{\Delta}
$
holds for all $i$ from $0$ through $N$, where
   $E_{k_i}$ is the unity $(k_i \times k_i)\,$-matrix
and
   $\mathit{\Delta} = \partial_1^2 + \ldots + \partial_n^2$.
Let $e (x)$ be the standard fundamental solution of convolution type for the
Helmholtz operator $\mathit{\Delta} + k^2$ in $\mathbb{R}^n$, that vanishes at
$\infty$.
E.g., if $n = 3$ then
$$
   e (x)
 = \frac{-1}{4 \pi}\,
   \frac{\exp (\imath k |x|)}{|x|}
$$
away from the origin.

\begin{corollary}
\label{c.SC}
Let
   (\ref{eq.complex}) be a complex of differential operators with constant
   coefficients of Dirac type,
and
   $\mathcal{X}$ a closed bounded domain with smooth boundary
in $\mathbb{R}^n$.
If
   $u \in H^1 (\mathcal{X}, F^i)$ and
   $f \in H^1 (\mathcal{X}, F^{i+1})$
satisfy $M^i (u,f) = 0$ in the interior of $\mathcal{X}$, then
\begin{equation}
\label{eq.SC}
   \Big(
         \begin{array}{cc}
           (1/\imath k) A^i{}^\ast A^i
         & - A^i{}^\ast
\\
           - A^i
         & - (1/\imath k) A^i A^i{}^\ast
         \end{array}
   \Big)
   \int_{\partial \mathcal{X}}
   e (x-\cdot)
   \binom {\imath n (f)}
          {\imath \sigma (A^i) (\nu) t (u)}
   ds
 = \binom {u (x)}{f (x)}
\end{equation}
   for all $x \in \mathcal{X} \setminus \partial \mathcal{X}$,
and the left-hand side vanishes away from $\mathcal{X}$.
\end{corollary}

Formula (\ref{eq.SC}) is in terms of the only operator $A^i$ which enters into
the definition of Maxwell's equations.
It extends modulo smoothing operators to general Maxwell's equations.

\begin{proof}
For a Dirac type complex one can choose
   $G^i u = - e \ast u$
for distribution sections $u$ of $F^i$ with compact support in $\mathbb{R}^n$,
where by $e \ast u$ is meant the convolution of $e$ and $u$.
The Schwartz kernel of $G^i$ is obviously $E_{k_i} e (x-y)$.
Hence, the Schwartz kernel of $\mathit{\Phi}^i$ is
$$
   - C^i{}' (y,\partial_y) \Big( E_{k_i + k_{i+1}} e (x-y) \Big)
 = - C^i (y,\partial_x) e (x-y),
$$
the variable $y$ of $C^i$ can be neglected, for the coefficients are constant.
From (\ref{eq.quotient}) one easily deduces that
$$
   C^i
 = \Big(
         \begin{array}{cc}
           (1/\imath k) (\mathit{\Delta}^i - k^2 - A^i{}^\ast A^i)

         & A^i{}^\ast
\\
           A^i
         & - (1/\imath k) (\mathit{\Delta}^{i+1} - k^2 - A^i A^i{}^\ast)
         \end{array}
   \Big),
$$
and so
$$
   \mathit{\Phi}^i (\cdot,y)
 = -
   \Big(
         \begin{array}{cc}
           - (1/\imath k) A^i{}^\ast A^i e (\cdot-y)

         & A^i{}^\ast e (\cdot-y)
\\
           A^i e (\cdot-y)
         & (1/\imath k) A^i A^i{}^\ast e (\cdot-y)
         \end{array}
   \Big)
$$
on $\mathbb{R}^n \setminus \{ y \}$, for
   $(\mathit{\Delta} + k^2) e = 0$
away from the origin.
To complete the proof it suffices to use Theorem \ref{t.Green} and observe
that
$
   \sigma (\mathit{\Delta}^i) (\nu)
 = |\nu|^2 E_{k_i}
 = E_{k_i}
$
in our particular case.
\end{proof}

As is easy to check,
$$
   M^i
   \left(
         \begin{array}{cc}
           \displaystyle
           \frac{1}{\imath k} A^i{}^\ast A^i
         & - A^i{}^\ast
\\
           - A^i
         & \displaystyle
           - \frac{1}{\imath k} A^i A^i{}^\ast
         \end{array}
   \right)
 =
   \left(
         \begin{array}{cc}
           0
         & \displaystyle
           - \frac{1}{\imath k} A^i{}^\ast (\mathit{\Delta}^{i+1} - k^2)
\\
           \displaystyle
           \frac{1}{\imath k} A^i (\mathit{\Delta}^i - k^2)
         & 0
         \end{array}
   \right),
$$
hence the left-hand side of (\ref{eq.SC}) satisfies Maxwell's equations in
   $\mathbb{R}^n \setminus \partial \mathcal{X}$
for all integrable functions $t (u)$ and
                             $n (f)$
on the boundary.
In this way we arrive at what is usually referred to as
   the Cauchy type integral
associated to Maxwell's equations,
   cf. \cite[3.2.3]{Tark90}.

\section{Boundary value problems}
\label{s.bvp}

The Stratton-Chu formula of Theorem \ref{t.Green} manifests the Cauchy data
on $\partial \mathcal{X}$ of sections $u \in H^1 (\mathcal{X},F^i)$ and
                                      $f \in H^1 (\mathcal{X},F^{i+1})$
with respect to Maxwell's operator $M^i$.
These are $t (u)$ and
          $n (f)$.

\begin{lemma}
\label{l.uniqueness}
Suppose $k^2$ is not real.
Let $u \in C^{2} (\mathcal{X},F^i)$ and
    $f \in C^{2} (\mathcal{X},F^{i+1})$
satisfy $M^i (u,f) = 0$ in $\mathcal{X}$.
Either of the conditions
   $t (u) = 0$ and
   $n (f) = 0$
on $\partial \mathcal{X}$ implies $(u,f) \equiv 0$ in all of $\mathcal{X}$.
\end{lemma}

\begin{proof}
Suppose $t (u) = 0$ on $\partial \mathcal{X}$.
Then $(Au,Au) = (A^\ast A u,u)$.
The second Maxwell equation yields
   $A^\ast A u = (\imath k) A^\ast f$
in $\mathcal{X}$.
By the first Maxwell equation,
   $(\imath k) A^\ast f = k^2 u$
in $\mathcal{X}$.
Hence $(Au,Au) = k^2 (u,u)$.
Since the imaginary part of $k^2$ does not vanish, we get $u \equiv 0$ in
$\mathcal{X}$.
It follows that $f \equiv 0$ in $\mathcal{X}$, showing the first part of the
lemma.
The proof of the second part is similar because of the obvious symmetry.
\end{proof}

Both $t (u)$ and
     $n (f)$
are sections of a vector bundle $F^i_t$ over the boundary.
This is the kernel of the bundle homomorphism
   $\sigma (A^{i-1})^\ast (\cdot,\nu)$
acting from $F^i |_{\partial \mathcal{X}}$
         to $F^{i-1} |_{\partial \mathcal{X}}$,
cf. \cite[3.2.2]{Tark90}.

It is not our purpose to study general boundary value problems for solutions
of Maxwell's equations.
We merely discuss a particular problem which stems from scattering of incident
electromagnetic waves by a perfectly conducting body.
In this case the tangential component $t (u)$ of the ``electric'' field $u$
must vanish on the body surface $\partial \mathcal{X}$.

Given a section $u_0$ of $F^i_t$ over $\partial \mathcal{X}$, we consider the
problem of finding a solution $(u,f)$ of Maxwell's equations in the interior
of $\mathcal{X}$, such that the tangential part of $u$ on the boundary is well
defined and $t (u) = u_0$ on $\partial \mathcal{X}$.
If a solution exists and is sufficiently smooth in $\mathcal{X}$, then
   $A^\ast u = 0$
and so
   $\mathit{\Delta}^i u = A^\ast A u = k^2 u$
in $\mathcal{X}$.
The Dirichlet problem for solutions of the elliptic equation
   $(\mathit{\Delta}^i - k^2) u = 0$
in $\mathcal{X}$ requires not only given values
   $t (u)$ on the boundary $\partial \mathcal{X}$
but also
   $n (u)$.
The homogeneous formal adjoint problem fails to have nonzero solutions unless
$\bar{k}{}^2$ (and so $k$) belongs to a discrete set of nonnegative numbers
with the only accumulation point at infinity.
In this latter case the Dirichlet problem is solvable if and only if the pair
   $(u_0,n (u))$
satisfies some orthogonality conditions.
In any case the canonical solution to the Dirichlet problem is written in the
form
   $u = \wp^i (u_0, n (u))$,
where $\wp^i$ is the so-called Poisson operator.
On having granted the component $u$ in $\mathcal{X}$ we define the component
$f$ from the second Maxwell equation by
   $f = (1/\imath k) Au$.
In order that the first Maxwell equation be satisfied it is necessary and
                                                            sufficient
that $A^\ast u = 0$, for then
\begin{eqnarray*}
   \imath k\, u + A^\ast f
 & = &
   \imath k\, u + (1/\imath k) A^\ast Au
\\
 & = &
   (1/\imath k) (\mathit{\Delta}^i - k^2) u
\\
 & = &
   0
\end{eqnarray*}
in $\mathcal{X}$.
Summarising we conclude that the above boundary value problem is solvable if
and only if
   $A^{\ast} \wp^i (u_0, n (u)) = 0$
holds in $\mathcal{X}$.
Since $\mathit{\Delta}$ and
      $A^{\ast}$
commute and so
   $(\mathit{\Delta}^{i-1} - k^2) A^{\ast} u = 0$
in the interior of $\mathcal{X}$, the equality
   $A^{\ast} \wp^i (u_0, n (u)) = 0$
is fulfilled in $\mathcal{X}$ if and only if it is fulfilled on the boundary.
In this way we arrive at the following result.

\begin{theorem}
\label{t.bvp}
As formulated above, the boundary value problem is equivalent to the equation
   $A^{\ast} \wp^i (u_0, n (u)) = 0$
on $\partial \mathcal{X}$.
\end{theorem}

For a treatment of this boundary integral equation in electromagnetic
scattering theory we refer to \cite{PikeSaba02}.

To demonstrate Theorem \ref{t.bvp} we look more closely at the extreme cases
   $i = 0$ and
   $i = N$.
For $i = 0$, the operator $A^{i-1}$ is zero, and so the normal part $n (u)$
vanishes and the equation
   $A^{\ast} \wp^i (u_0, n (u)) = 0$
on the boundary is automatically fulfilled.
The theorem says that the boundary value problem in question has a solution
for each
   $u_0 \in \mathcal{D}' (\partial \mathcal{X}, F^0)$,
and the canonical solution is actually given by $u = \wp^0 (u_0)$.
If $i = N$, then the operator $A^{i-1}{}^\ast$ has injective symbol whence the
bundle $F^i_t$ is zero.
Thus, the boundary condition $t (u) = u_0$ is empty.
However, the section $F^{N+1}$ is zero and so $f \equiv 0$, which implies
                                              $u \equiv 0$
in $\mathcal{X}$.
This corresponds to the fact that the boundary integral equation
   $A^{\ast} \wp^i (u) = 0$
for a section $u$ of $F^N |_{\partial \mathcal{X}}$ has the only solution
   $u = 0$.

\section{The Cauchy problem}
\label{s.tCp}

When combined with Maxwell's equations, the data $t (u)$ and
                                                 $n (f)$
actually determine both $u$ and
                       $f$
on the boundary, for
$$
   \begin{array}{rcl}
     n (u)
   & =
   & - (1/\imath k) n (A^\ast f),
\\
     t (f)
   & =
   & (1/\imath k) t (Au)
   \end{array}
$$
while $n (A^\ast f)$ and
      $t (Au)$
are uniquely determined by $n (f)$ and
                           $t (u)$
(the so-called tangential operators on $\partial \mathcal{X}$, cf.
   \cite[3.1.5]{Tark90}).
Hence, the Cauchy problem for Maxwell's equations in $\mathcal{X}$ consists in
finding a solution $(u,f)$ to $M^i (u,f) = 0$ with prescribed data $t (u)$ and
                                                                   $n (u)$
on a part of the surface $\partial \mathcal{X}$.
For a general theory the reader can consult
   \cite{Tark95}.
To obtain constructive results we restrict ourselves to Dirac type complexes
of differential operators with constant coefficients in $\mathbb{R}^n$, as
defined in Section \ref{s.SCf}.

Let $\mathcal{S}$ be an open piece on the boundary of $\mathcal{X}$.
The Cauchy problem for solutions of Maxwell's equations in $\mathcal{X}$ with
data on $\mathcal{S}$ consists in the following.
Given sections $u_0$ and
               $f_0$
of the bundle $F^i_t$ over $\mathcal{S}$, find a solution $(u,f)$ to
   $M^i (u,f) = 0$
in the interior of $\mathcal{X}$, such that
   $t (u) = u_0$ and
   $n (u) = f_0$
on $\mathcal{S}$.
To study this problem we have to choose function spaces for $u_0$,
                                                            $f_0$
and $(u,f)$.
It is a little cumbersome, for the behaviour of the solution near the boundary
of $\mathcal{S}$ requires a careful study.
In order to highlight principal difficulties in the Cauchy problem we restrict
our attention to the case where $u_0$ and
                                $f_0$
are continuous sections of $F^i_t$ over the closure of $\mathcal{S}$.
If exists, the solution $(u,f)$ should be continuous up to $\mathcal{S}$ in
the interior of $\mathcal{X}$, and we avoid discussion of weak boundary values
of $(u,f)$ on $\mathcal{S}$.
We thus consider the problem
\begin{equation}
\label{eq.tCp}
\left\{ \begin{array}{rclcl}
          M^i (u,f)
        & =
        & 0
        & \mbox{in} & \mathcal{X} \setminus \partial \mathcal{X},
\\
          t (u)
        & =
        & u_0
        & \mbox{on} & \mathcal{S},
\\
          n (f)
        & =
        & f_0
        & \mbox{on} & \mathcal{S}.
        \end{array}
\right.
\end{equation}

It is well known that this problem has at most one solution in any reasonable
space of functions on $\mathcal{X}$.

To study (\ref{eq.tCp}) we introduce an integral completely determined by the
Cauchy data $u_0$ and
            $f_0$
on $\mathcal{S}$,
   namely
\begin{equation}
\label{eq.dlp}
   \mathcal{G} (u_0,f_0) (x)
 = \int_{\mathcal{S}}
   \mathit{\Phi}^i (x-\cdot)
   \binom {\imath f_0}
          {\imath \sigma (A^i) (\nu) u_0}
   ds
\end{equation}
for $x \in \mathbb{R}^n \setminus \overline{\mathcal{S}}$.
Since the fundamental solution $\mathit{\Phi}^i (x-y)$ is real analytic away
from the diagonal of $\mathbb{R}^n \times \mathbb{R}^n$, it follows that
   $\mathcal{G} (u_0,f_0)$
is real analytic in the complement of $\overline{\mathcal{S}}$.

Moreover, the potential (\ref{eq.dlp}) satisfies Maxwell's equations
   $M^i \mathcal{G} (u_0,f_0) = 0$
in $\mathbb{R}^n \setminus \overline{\mathcal{S}}$.
In particular, the components of the vector-valued function
   $\mathcal{G} (u_0,f_0)$
prove to be solutions to the (scalar) Helmholtz equation
   $(\mathit{\Delta} + k^2) \mathcal{G} = 0$
in the complement of $\overline{\mathcal{S}}$.

When $x$ crosses the hypersurface $\mathcal{S}$, the integral
   $\mathcal{G} (u_0,f_0) (x)$
has jump.
The corresponding jump formulas look very like
   the classical Sokhotskii-Plemelj formulas,
cf.
   \cite[3.2.3]{Tark90}.

\begin{theorem}
\label{t.scftCp}
In order that there be a solution $(u,f)$ of the Cauchy problem (\ref{eq.tCp})
continuous up to $\mathcal{S}$, it is necessary and sufficient that the
integral
   $\mathcal{G} (u_0,f_0)$
might be extended from $\mathbb{R}^n \setminus \mathcal{X}$
                  through $\mathcal{S}$
                  to the interior of $\mathcal{X}$
as a real analytic function.
\end{theorem}

\begin{proof}
{\it Necessity.}
Suppose there is a solution $(u,f)$ of the Cauchy problem (\ref{eq.tCp}) which
is continuous up to $\mathcal{S}$.
Define sections $U$ and
                $F$
of $F^i$ and
   $F^{i+1}$
in $\mathbb{R}^n \setminus \partial \mathcal{X}$ by
$$
   \binom {U}{F}
 = \left\{ \begin{array}{lcl}
             \displaystyle
             \mathcal{G} (u_0,f_0) - \binom {u}{f}
           & \mbox{in}
           & \mathcal{X} \setminus \partial \mathcal{X},
\\
             \mathcal{G} (u_0,f_0)
           & \mbox{in}
           & \mathbb{R}^n \setminus \mathcal{X}.
           \end{array}
   \right.
$$
Write $U^{\pm}$ for the restriction of $U$ to the open sets
   $\mathcal{X} \setminus \partial \mathcal{X}$ and
   $\mathbb{R}^n \setminus \mathcal{X}$,
respectively, and similarly for $F$.

We can assume, by shrinking $\mathcal{X}$ if necessary, that both
   $u$ and
   $f$
are continuous up to the boundary of $\mathcal{X}$.
Using Theorem \ref{t.Green}, we get
$$
   \binom {U^{+} (x)}{F^{+} (x)}
 =
 - \int_{\partial \mathcal{X} \setminus \mathcal{S}}
   \mathit{\Phi}^i (x-\cdot)
   \binom {\imath n (f)}
          {\imath \sigma (A^i) (\nu) t (u)}
   ds
$$
for all $x$ in the interior of $\mathcal{X}$.
Hence it follows that $U^{+}$ and
                      $F^{+}$
extend through $\mathcal{S}$ to real analytic sections of $F^i$ and
                                                          $F^{i+1}$
on all of
   $(\mathbb{R}^n \setminus \partial \mathcal{X}) \cup \mathcal{S}$.
Furthermore, applying Theorem \ref{t.Green} once again yields
\begin{eqnarray*}
 - \int_{\partial \mathcal{X} \setminus \mathcal{S}}
   \mathit{\Phi}^i (x,\cdot)
   \binom {\imath n (f)}
          {\imath \sigma (A^i) (\nu) t (u)}
   ds
 & = &
   \int_{\mathcal{S}}
   \mathit{\Phi}^i (x,\cdot)
   \binom {\imath n (f)}
          {\imath \sigma (A^i) (\nu) t (u)}
   ds
\\
 & = &
   \binom {U^{-} (x)}{F^{-} (x)}
\end{eqnarray*}
for
   $x \in \mathbb{R}^n \setminus \mathcal{X}$.
Therefore,
   $\mathcal{G} (u_0,f_0)$
extends from $\mathbb{R}^n \setminus \mathcal{X}$
        through $\mathcal{S}$
        to the interior of $\mathcal{X}$
as a real analytic function, as desired.

{\it Sufficiency.}
Conversely, let $U$ and
                $F$
be real analytic sections of bundles $F^i$ and
                                     $F^{i+1}$
over
   $(\mathbb{R}^n \setminus \partial \mathcal{X}) \cup \mathcal{S}$,
such that
$$
   \binom {U}{F} = \mathcal{G} (u_0,f_0)
$$
away from $\mathcal{X}$.
Then $M^i (U,F) = 0$ in $\mathbb{R}^n \setminus \mathcal{X}$.
Since $M^i (U,F)$ is real analytic in
   $(\mathbb{R}^n \setminus \partial \mathcal{X}) \cup \mathcal{S}$,
it actually vanishes in $\mathcal{X}$, too.

Set
$$
   \binom {u (x)}{f (x)}
 = \mathcal{G} (u_0,f_0) (x) - \binom {U (x)}{F (x)}
$$
for $x$ in the interior of $\mathcal{X}$.
From what has already been proved it follows that
   $u$ and
   $f$
are continuous up to $\mathcal{S}$ and satisfy Maxwell's equations.
We claim that $(u,f)$ is the desired solution of (\ref{eq.tCp}).
To see this, it remains to verify that
   $t (u) = u_0$ and
   $n (f) = f_0$
on $\mathcal{S}$.

To this end, denote by $\mathcal{G}_e (u_0,f_0)$ and
                       $\mathcal{G}_m (u_0,f_0)$
the components of $\mathcal{G} (u_0,f_0)$ that are sections of the bundles
   $F^i$ and
   $F^{i+1}$,
respectively (electric and magnetic fields).
Since $U$ and
      $F$
are continuous in
   $(\mathbb{R}^n \setminus \partial \mathcal{X}) \cup \mathcal{S}$,
we get by the Sokhotskii-Plemelj formulas
\begin{eqnarray*}
   t (u)
 & = &
   t \left( \mathcal{G}_e (u_0,f_0)^{+} - U^{+} \right)
\\
 & = &
   t \left( \mathcal{G}_e (u_0,f_0)^{+} - U^{-} \right)
\\
 & = &
   t \left( \mathcal{G}_e (u_0,f_0)^{+} - \mathcal{G}_e (u_0,f_0)^{-} \right)
\\
 & = &
   u_0
\end{eqnarray*}
on $\mathcal{S}$.
Analogously,
\begin{eqnarray*}
   n (f)
 & = &
   n \left( \mathcal{G}_m (u_0,f_0)^{+} - F^{+} \right)
\\
 & = &
   n \left( \mathcal{G}_m (u_0,f_0)^{+} - F^{-} \right)
\\
 & = &
   n \left( \mathcal{G}_m (u_0,f_0)^{+} - \mathcal{G}_m (u_0,f_0)^{-} \right)
\\
 & = &
   f_0
\end{eqnarray*}
on $\mathcal{S}$, which is our assertion.
\end{proof}

Theorem \ref{t.scftCp} may be summarised by saying that any explicit formula
for analytic continuation from $\mathbb{R}^n \setminus \mathcal{X}$
                          through $\mathcal{S}$
                          to the interior of $\mathcal{X}$
leads to a formula for solutions of the Cauchy problem with data on
$\mathcal{S}$.
This idea goes back at least as far as
   \cite{Aize93}.

\section{Expansion of the fundamental solution}
\label{s.Eotfs}

In this section we extend the results of Section 6 of \cite{MakhNiyoTark08}
from the case $n = 3$ to the case of arbitrary $n > 3$.

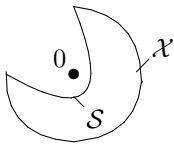
\begin{figure}[h]
\special{em:linewidth 0.4pt}
\unitlength 1mm
\linethickness{0.4pt}
\begin{picture}(20.33,19.00)
\put(10.00,10.00){\circle*{1.33}}
\emline{10.00}{19.00}{1}{10.84}{16.76}{2}
\emline{10.84}{16.76}{3}{11.49}{14.75}{4}
\emline{11.49}{14.75}{5}{11.95}{12.97}{6}
\emline{11.95}{12.97}{7}{12.20}{11.42}{8}
\emline{12.20}{11.42}{9}{12.27}{10.09}{10}
\emline{12.27}{10.09}{11}{12.14}{8.99}{12}
\emline{12.14}{8.99}{13}{11.82}{8.12}{14}
\emline{11.82}{8.12}{15}{11.30}{7.47}{16}
\emline{11.30}{7.47}{17}{10.58}{7.05}{18}
\emline{10.58}{7.05}{19}{9.68}{6.86}{20}
\emline{9.68}{6.86}{21}{8.57}{6.90}{22}
\emline{8.57}{6.90}{23}{7.28}{7.17}{24}
\emline{7.28}{7.17}{25}{5.79}{7.66}{26}
\emline{5.79}{7.66}{27}{4.10}{8.38}{28}
\emline{4.10}{8.38}{29}{1.00}{10.00}{30}
\emline{10.00}{19.00}{31}{12.33}{18.57}{32}
\emline{12.33}{18.57}{33}{14.31}{17.83}{34}
\emline{14.31}{17.83}{35}{15.94}{16.78}{36}
\emline{15.94}{16.78}{37}{17.22}{15.42}{38}
\emline{17.22}{15.42}{39}{18.16}{13.75}{40}
\emline{18.16}{13.75}{41}{18.75}{11.78}{42}
\emline{18.75}{11.78}{43}{19.00}{9.00}{44}
\emline{19.00}{9.00}{45}{19.00}{9.00}{46}
\emline{19.00}{9.00}{47}{19.00}{9.00}{48}
\emline{19.00}{9.00}{49}{18.60}{6.87}{50}
\emline{18.60}{6.87}{51}{17.83}{5.07}{52}
\emline{17.83}{5.07}{53}{16.69}{3.59}{54}
\emline{16.69}{3.59}{55}{15.17}{2.44}{56}
\emline{15.17}{2.44}{57}{13.28}{1.61}{58}
\emline{13.28}{1.61}{59}{10.00}{1.00}{60}
\emline{10.00}{1.00}{61}{7.59}{1.34}{62}
\emline{7.59}{1.34}{63}{5.56}{1.99}{64}
\emline{5.56}{1.99}{65}{3.90}{2.94}{66}
\emline{3.90}{2.94}{67}{2.61}{4.21}{68}
\emline{2.61}{4.21}{69}{1.70}{5.78}{70}
\emline{1.70}{5.78}{71}{1.16}{7.67}{72}
\emline{1.16}{7.67}{73}{1.00}{10.33}{74}
\put(9.00,11.00){\makebox(0,0)[rb]{$0$}}
\emline{10.00}{7.00}{75}{11.33}{5.67}{76}
\put(11.67,5.33){\makebox(0,0)[lt]{$\mathcal{S}$}}
\emline{18.33}{10.33}{77}{20.00}{12.00}{78}
\put(20.33,12.33){\makebox(0,0)[lb]{$\mathcal{X}$}}
\end{picture}
\caption{A typical domain under consideration.}
\label{f.lune}
\end{figure}

We apply the method of bases with double orthogonality to the Cauchy problem
(\ref{eq.tCp}) in the particular case, where $\mathcal{X}$ is a part of a ball
$B = B (0,R)$ with centre at the origin and radius $R > 0$.
Let $\mathcal{S}$ be a smooth closed hypersurface in $B$ which does not meet
$x = 0$ and divides $B$ into two domains.
Denote by $\mathcal{X}$ the closed domain that does not contain the origin.
Its boundary $\partial \mathcal{X}$ consists of $\mathcal{S}$ and a part of
the sphere $\partial B$ in $\mathbb{R}^{n}$,
   cf. Figure~\ref{f.lune}.
The advantage of using domains $\mathcal{X}$ of the above form lies in the
fact that the problem reduces to analytic continuation
      from a small ball around $0$
      to $B$.

By the above, the double layer potential $\mathcal{G} (u_0,f_0)$ satisfies
both
   $M^i \mathcal{G} (u_0,f_0) = 0$
and
   $(\mathit{\Delta} + k^2) \mathcal{G} (u_0,f_0) = 0$
away from the closure of $\mathcal{S}$.
The latter equation is the scalar Helmholtz equation and follows from the
former one.
The idea is now to use bases with double orthogonality in Hilbert spaces of
solutions to the Helmholtz equation to derive conditions of analytic
continuation
   from $B (0,\epsilon)$
   to $B (0,R)$
for $\mathcal{G} (u_0,f_0)$.

In spherical coordinates, the Helmholtz operator $\mathit{\Delta} + k^2$ in
the space $\mathbb{R}^{n}$ takes the form
\begin{equation}
\label{eq.Hoispc}
   \mathit{\Delta} + k^2
 = \frac{1}{r^{2}}
   \Big( \Big( r \frac{\partial}{\partial r} \Big)^2
       + (n-2)\, r \frac{\partial}{\partial r}
       + k^{2} r^{2}
       - \mathit{\Delta}_{\mathbb{S}^{n-1}}
   \Big),
\end{equation}
where
   $\mathit{\Delta}_{\mathbb{S}^{n-1}}$ is the Laplace-Beltrami operator on
   the unit sphere.
Recall that $k$ is an arbitrary complex number with $\Im k \geq 0$.

To solve the homogeneous equation $(\mathit{\Delta} + k^2) u = 0$ we use the
Fourier method of separation of variables.
Writing
   $u (r,\omega) = g (r,k) h (\omega)$
we get two separate equations for $g$ and
                                  $h$,
namely
\begin{eqnarray*}
   \Big( \Big( r \frac{\partial}{\partial r} \Big)^2
       + (n-2)\, r \frac{\partial}{\partial r}
       + k^{2} r^{2}
   \Big) g
 & = &
   c\, g
\\
   \Delta_{\mathbb{S}^{n-1}} h
 & = &
   c\, h,
\end{eqnarray*}
$c$ being an arbitrary constant.

The second equation has nonzero solutions if and only if $c$ is an eigenvalue
of $\mathit{\Delta}_{\mathbb{S}^{n-1}}$.
These are well known to be
   $c = \nu (\nu + n - 2)$,
for $\nu = 0, 1, \ldots$, cf.
   \cite{TikhSama72} and elsewhere.
The corresponding eigenfunctions of $\mathit{\Delta}_{\mathbb{S}^{n-1}}$ are
spherical harmonics $h_\nu (\omega)$ of degree $\nu$, i.e.
\begin{equation}
\label{eq.spherharm}
   \mathit{\Delta}_{\mathbb{S}^{n-1}} h_\nu
 = \nu (\nu + n - 2)\, h_\nu.
\end{equation}

Consider now the ordinary differential equation with respect to the variable
$r > 0$
\begin{equation}
\label{eq.Bessel}
   \Big( \Big( r \frac{\partial}{\partial r} \Big)^2
       + (n-2)\, r \frac{\partial}{\partial r}
       + \left( k^2 r^2 - \nu (\nu + n - 2) \right)
   \Big) g (r,k)
 = 0.
\end{equation}
This is a version of the Bessel equation, and the space of its solutions is
two-dimensional.

For example, if $n = 3$ and
                $k = 0$,
then
   $g (r,0) = a r^\nu + b r^{-\nu-1}$
with arbitrary constants $a$ and
                         $b$
is a general solution to (\ref{eq.Bessel}).
In this situation any function $r^\nu h_\nu (\varphi)$ is a homogeneous
harmonic polynomial.
In the general case the space of solutions to (\ref{eq.Bessel}) contains
a one-dimensional subspace of functions bounded at the point $r = 0$, cf.
   \cite{TikhSama72}.

For $\nu = 0, 1, \ldots$, fix a non-zero solution $g_\nu (r,k)$ of
(\ref{eq.Bessel}) which is bounded at $r = 0$.
Then
\begin{equation}
\label{eq.Helmholtz}
   \left( \mathit{\Delta} + k^2 \right)
   \left( g_\nu (r,k) h_\nu (\omega) \right)
 = 0
\end{equation}
on all of $\mathbb{R}^n$.
Indeed, by (\ref{eq.Hoispc}),
           (\ref{eq.spherharm}) and
           (\ref{eq.Bessel})
we conclude that this equality holds in $\mathbb{R}^n \setminus \{ 0 \}$.
We now use the fact that $g_\nu (r,k) h_\nu (\omega)$ is bounded at the origin
to see that (\ref{eq.Helmholtz}) holds.

It is known that there are exactly $J (\nu)$ linearly independent spherical
harmonics of degree $\nu$, where
$$
   J (\nu)
 = \frac{(2\nu\!+\!n\!-\!2) (\nu\!+\!n\!-\!3) !}{(n\!-\!2) ! \nu !}.
$$
Pick an orthonormal basis
$$
   \{ h_\nu^{(j)} \}_{\nu = 0, 1, \ldots \atop j = 1, \ldots, J (\nu)}
$$
in $L^{2} (\mathbb{S}^{n-1})$.

\begin{lemma}
\label{l.basis}
For every $R > 0$, the system
\begin{equation}
\label{eq.bij}
   \{ b_{\nu}^{(j)} (r, \omega, k) := g_\nu (r,k) h_\nu^{(j)} (\omega)
   \}_{\nu = 0, 1, \ldots \atop j = 1, \ldots, J (\nu)}
\end{equation}
is an orthogonal basis in the subspace of $L^{2} (B (0,R))$ consisting of
solutions to the Helmholtz equation $(\mathit{\Delta} + k^{2}) u = 0$.
\end{lemma}

\begin{proof}
Indeed, as
   $\{ h_{\nu}^{(j)} \}$
is an orthonormal basis in the space $L^2 (\mathbb{S}^{n-1})$ on the unit
sphere, the system (\ref{eq.bij}) is orthogonal in $L^2 (B (0,R))$ because
\begin{eqnarray*}
   (b_{\mu}^{(i)}, b_{\nu}^{(j)})_{L^2 (B (0,R))}
 & = &
   (h_{\mu}^{(i)}, h_{\nu}^{(j)})_{L^2 (\mathbb{S}^{n-1})}
   \int_0^R g_\mu (r,k) \overline{g_\nu (r,k)}\, r^{n-1} dr                                                 \\
 & = &
   0
\end{eqnarray*}
for $\mu \ne \nu$ or
    $i \ne j$.
Finally, since the system of harmonics $\{ h_\nu^{(j)} \}$ is dense in
                                       $C^{\infty} (\mathbb{S}^{n-1})$
we see that system (\ref{eq.bij}) is dense in the subspace of
   $L^{2} (B (0,R))$
consisting of the solutions of the Helmholtz equation in the ball,
   which completes the proof.
\end{proof}

For any fixed $y \in \mathbb{R}^n$ with $|y| > R$, the fundamental solution
$e (x-y)$ of the Helmholtz equation satisfies
   $(\mathit{\Delta} + k^2) e (\cdot-y) = 0$
in $B (0,R)$ and is obviously square integrable in the ball.
Therefore, $e (\cdot-y)$ can be represented in $B (0,R)$ by its Fourier series
with respect to the orthogonal system (\ref{eq.bij})
\begin{equation}
\label{eq.eofs}
   e (x-y)
 = \sum_{\nu=0}^{\infty}
   \sum_{j=1}^{J (\nu)}
   c_\nu^{(j)} (y,k)\, b_\nu^{(j)} (x,k),
\end{equation}
the series converges in the $L^2 (B (0,R))\,$-norm.
The Fourier coefficients $c_\nu^{(j)} (y,k)$ are defined by familiar formulas
through the scalar product in $L^2 (B (0,R))$.
These formulas show that $c_\nu^{(j)} (y,k)$ are real analytic functions of
$y$ provided $|y| > R$.
Moreover, they satisfy the Helmholtz equation in the complement of
   $\overline{B (0,R)}$.
On the other hand, the system (\ref{eq.bij}) is orthogonal in any space
   $L^{2} (B (0,R))$
with $R > 0$.
Hence, the coefficients $c_\nu^{(j)} (y,k)$ are actually independent of the
particular choice of $R$ satisfying $R < |y|$.
This shows that $c_\nu^{(j)} (y,k)$ are solutions of the Helmholtz equation
away from the origin in $\mathbb{R}^n$.
Since $e (x-y) = O (|x-y|^{2-n})$ as $x \to y$, we may tend $R \to |y|$ in the
formulas for $c_\nu^{(j)} (y,k)$.
This immediately yields explicit equalities
$$
   c_\nu^{(j)} (y,k)
 = \Big( e (\cdot-y), b_\nu^{(j)} (\cdot,k) \Big)_{L^{2} (B (0,|y|))}
   \, \Big/ \,
   \int_0^{|y|} |g_\nu (r,k)|^{2} r^{n-1} dr
$$
for $\nu = 0, 1, \ldots$ and
    $j = 1, \ldots, J (\nu)$.

\begin{lemma}
\label{l.eofs}
In the cone
   $\{ (x,y) \in \mathbb{R}^{n} \times \mathbb{R}^{n} : |x| / |y| < 1 \}$,
one has the Fourier series expansion (\ref{eq.eofs}), where
   the series converges uniformly together with all its derivatives on compact
   subsets of the cone.
\end{lemma}

\begin{proof}
The Fourier series expansion is a direct consequence of
   Lemma \ref{l.basis}.
The uniform convergence on compact subsets of the cone presents a more
delicate problem.
It can be handled in a familiar way, for $b_\nu^{(j)}$ are solutions of the
Helmholtz equation.
\end{proof}

\section{Regularisation}
\label{s.reg}

Substituting (\ref{eq.eofs}) into the formula for the fundamental solution
$\mathit{\Phi}^i (x-y)$ of $M^i$ we obtain
\begin{equation}
\label{eq.eofsoM}
   \mathit{\Phi}^i (x-y)
 = \sum_{\nu=0}^{\infty}
   \mathit{\Phi}^i_{\nu} (x,y)
\end{equation}
where the series converges uniformly along with all its derivatives on compact
subsets of the cone
   $\{ (x,y) \in \mathbb{R}^{n} \times \mathbb{R}^{n} : |x| / |y| < 1 \}$,
and
$$
   \mathit{\Phi}^i_{\nu} (x,y)
 =
   \Big( \!
         \begin{array}{cc}
           (1/\imath k) A^i{}^\ast (\partial_x) A^i (\partial_x)
      \! & \!
           - A^i{}^\ast (\partial_x)
\\
           - A^i (\partial_x)
      \! & \!
           - (1/\imath k) A^i (\partial_x) A^i{}^\ast (\partial_x)
         \end{array}
         \!
   \Big)
   \sum_{j=1}^{J (\nu)} c_\nu^{(j)} (y,k)\, b_\nu^{(j)} (x,k)
$$
for $\nu = 0, 1, \ldots$.
We obtain the same formulas if we substitute $- \partial_y$
                         for the derivatives $\partial_x$.

\begin{lemma}
\label{l.poiPjinu}
Every term
   $\mathit{\Phi}_{\nu} (x,y)$
is a real analytic matrix-valued function on the set
   $\mathbb{R}^n \times (\mathbb{R}^n \setminus \{ 0 \})$,
satisfying
$$
\begin{array}{rcl}
            M^i (\partial_x) \mathit{\Phi}^i_{\nu} (x,y) & = & 0,
\\
   M^i{}' (\partial_y) (\mathit{\Phi}_{\nu}^i (x,y))^{T} & = & 0.
\end{array}
$$
\end{lemma}

\begin{proof}
These properties are obvious by the very construction of degenerate kernels
   $\mathit{\Phi}^i_{\nu} (x,y)$.
The singularity at $y = 0$ is due to
$$
   \int_0^{|y|} |g_\nu (r,k)|^{2} r^{n-1} dr.
$$
\end{proof}

A series expansion like (\ref{eq.eofsoM}) with terms $\mathit{\Phi}^i_{\nu}$
satisfying the transposed equation
   $M' (\partial_y) (\mathit{\Phi}^i_{\nu} (x,y))^{T} = 0$
is already sufficient to derive an explicit formula for solutions of the
Cauchy problem (\ref{eq.tCp}).
Set
$$
   R^i_N (x,y)
 = \mathit{\Phi}^i (x-y) - \sum_{\nu=0}^{N} \mathit{\Phi}^i_{\nu} (x,y)
$$
for $(x,y) \in \mathbb{R}^n \times (\mathbb{R}^n \setminus \{ 0 \})$.
In this way we obtain what is referred to as a Carleman function of the Cauchy
problem, cf. \cite[10.4]{Tark95}.

\begin{theorem}
\label{t.Cf}
Suppose
   $(u,f)$ is a solution to Maxwell's equations in the interior of
   $\mathcal{X}$ continuous up to the closure of $\mathcal{S}$.
Then
$$
   \binom {u (x)}{f (x)}
 = \lim_{N \to \infty}
   \int_{\mathcal{S}}
   R^i_N (x,\cdot)
   \binom {\imath n (f)}
          {\imath \sigma (A^i) (\nu) t (u)}
   ds
$$
for all $x \in \mathcal{X} \setminus \partial \mathcal{X}$.
\end{theorem}

\begin{proof}
We can assume,
   by approximating $\mathcal{X}$ by smaller domains whose boundaries
   intersect the boundary of $\mathcal{X}$ only in $\mathcal{S}$,
that $(u,f)$ is smooth in $\mathcal{X}$.
Combining Corollary \ref{c.SC} with Stokes' formula and
                                    Lemma \ref{l.poiPjinu}
we conclude easily that
$$
   \binom {u (x)}{f (x)}
 = \int_{\partial \mathcal{X}}
   R^i_N  (x,\cdot)
   \binom {\imath n (f)}
          {\imath \sigma (A^i) (\nu) t (u)}
   ds
$$
for each fixed $x \in \mathcal{X} \setminus \partial \mathcal{X}$ and
    all $N = 0, 1, \ldots$.
Let $N \to \infty$.
Since the series (\ref{eq.eofsoM}) converges on $\partial B (0,R)$ uniformly
in $x$ on compact subsets of $B (0,R)$, it follows that the part of the
boundary integral over
   $\partial \mathcal{X} \setminus \mathcal{S}$
tends to zero.
This establishes the desired formula.
\end{proof}

Let $(u,f)$ be a solution of the Cauchy problem (\ref{eq.tCp}).
A straightforward computation shows that
\begin{equation}
\label{eq.auxiliary}
   \int_{\mathcal{S}}
   R^i_N (x,\cdot)
   \binom {\imath n (f)}
          {\imath \sigma (A^i) (\nu) t (u)}
   ds
 = \mathcal{G} (u_0,f_0) (x)
 - \binom {U_N (x)}
          {F_N (x)}
\end{equation}
for $x \not\in \overline{\mathcal{S}}$,
   where
$$
   \binom {U_N (x)}{F_N (x)}
 = \sum_{\nu=0}^{N}
   \int_{\mathcal{S}}
   \mathit{\Phi}^i_{\nu} (x,\cdot)
   \binom {\imath n (f)}
          {\imath \sigma (A^i) (\nu) t (u)}
   ds.
$$

Write $\epsilon > 0$ for the distance between $\mathcal{S}$ and the origin.
If $x \in B (0,\varepsilon)$ then the left-hand side of (\ref{eq.auxiliary})
tends to zero,
   for the series (\ref{eq.eofsoM}) converges uniformly on $\mathcal{S}$.
It follows that the series (i.e. a sequence of partial sums)
   $\{ (U_N,F_N) \}$
converges to
   $\mathcal{G} (u_0,f_0)$
uniformly together with its derivatives on compact subsets of the ball
   $B (0,\epsilon)$.

The last observation and Theorem \ref{t.scftCp} fit together to yield certain
conditions of solvability for the Cauchy problem.

\begin{corollary}
\label{c.scftCp}
If the series $\{ (U_N,F_N) \}$ converges uniformly on compact subsets of the
ball $B (0,R)$, then the Cauchy problem (\ref{eq.tCp}) is solvable.
\end{corollary}

\begin{proof}
Since the terms of the series $\{ (U_N,F_N) \}$ are component-wise solutions
of the Helmholtz equation, it follows by the Stieltijes-Vitali theorem that
its sum $(U,F) = \lim (U_N,F_N)$ satisfies component-wise the same equation in
$B (0,R)$.
Hence, $(U,F)$ is a real analytic section of $F^i \oplus F^{i+1}$ on $B (0,R)$.
As $(U,F)$ actually coincides with
   $\mathcal{G} (u_0,f_0)$
in the smaller ball $B (0,\epsilon)$, the solvability of the Cauchy problem
follows from Theorem \ref{t.scftCp}, as desired.
\end{proof}

In many interesting cases the uniform convergence of the series
   $\{ (U_N,F_N) \}$
is not only sufficient but also necessary for the solvability of
(\ref{eq.tCp}),
   cf. \cite[2.9]{Shla96}.

\section{Propagation of electromagnetic waves}
\label{s.poemw}

We now turn to classical Maxwell's equation in a three-dimensional space, this
latter can be a three-dimensional manifold $\mathcal{X}'$ as well.
If $\mathcal{X}'$ is compact and closed, the theory is especially instructive,
for the so-called Bohr-Sommerfeld radiation conditions are no longer needed.
To demonstrate our constructions along more classical lines, we consider the
case $\mathcal{X}' = \mathbb{R}^3$.
As mentioned in Section \ref{s.aMe}, classical Maxwell's equations have the
form
$$
\begin{array}{rcl}
   \imath k\, E + d^\ast H
 & =
 & 0,
\\
   - \imath k\, H + dE
 & =
 & 0,
\end{array}
$$
$E$ and
$H$ being functions in a closed domain $\mathcal{X} \subset \mathbb{R}^3$ with
values in $\mathbb{R}^3$.
If $E$ is suitably specified within $1\,$-forms and
   $H$ within $2\,$-forms,
both the exterior derivative $d$ and
     its formal adjoint $d^\ast$
can be identified with the operator $\mathrm{curl}$ on vector fields in
   $\mathbb{R}^3$.

Applying Corollary \ref{c.SC} to classical Maxwell's equations we obtain the
Stratton-Chu formula \cite{Stra41}.

\begin{theorem}
\label{t.cSC}
Suppose $(E,H)$ is an electromagnetic wave in $\mathcal{X}$ whose
   electric component $E$ and
   magnetic component $H$
are both continuous up to the boundary.
Then
$$
   \Big(
         \begin{array}{cc}
           (1/\imath k) d^\ast d
         & - d^\ast
\\
           - d
         & - (1/\imath k) d d^\ast
         \end{array}
   \Big)
   \int_{\partial \mathcal{X}}
   \frac{-1}{4 \pi} \frac{\exp (\imath k |x-\cdot|)}{|x-\cdot|}
   \binom {\imath n (H)}
          {- \nu \wedge t (E)}
   ds
 = \binom {E (x)}{H (x)}
$$
   for all $x \in \mathcal{X} \setminus \partial \mathcal{X}$,
and the left-hand side vanishes away from $\mathcal{X}$.
\end{theorem}

Let
   $\mathcal{X}$ be a domain in $\mathbb{R}^3$ and
   $\mathcal{S}$ an open piece on $\partial \mathcal{X}$
of the form displayed in Figure~\ref{f.lune}.
We consider the problem of finding an electric field $E$ and
                                   a magnetic field $H$
in $\mathcal{X}$ with given tangential component $E_0$ of $E$ and
                            normal component $H_0$ of $H$
on $\mathcal{S}$.

Since the fundamental solution of the Helmholtz equation in $\mathbb{R}^3$ is
square integrable, the Fourier series expansion (\ref{eq.eofs}) is evident.
The coefficients $c_\nu^{(j)} (y,k)$ are explicitly given by
$$
   c_\nu^{(j)} (y,k)
 = \Big( \frac{-1}{4 \pi} \frac{\exp (\imath k |\cdot-y|)}{|\cdot-y|},
         b_\nu^{(j)} (\cdot,k) \Big)_{L^{2} (B (0,|y|))}
   \, \Big/ \,
   \int_0^{|y|} |g_\nu (r,k)|^{2} r^{2} dr
$$
for $\nu = 0, 1, \ldots$ and
    $j = 1, \ldots, J (\nu)$,
where $J (\nu) = 2 \nu + 1$.

Now the theory of Section \ref{s.reg} gives a Carleman function $R_N (x,y)$ of
the Cauchy problem in explicit form.
Namely,
$$
   R_N (x,y)
 =
   \Big(
         \begin{array}{cc}
           (1/\imath k) d^\ast d
         & - d^\ast
\\
           - d
         & - (1/\imath k) d d^\ast
         \end{array}
   \Big)
   \Big(
   e (x-y)
 - \sum_{\nu=0}^{N}
   \sum_{j=1}^{J (\nu)}
   c_\nu^{(j)} (y,k)\, b_\nu^{(j)} (x,k)
   \Big),
$$
the differential operator on the right-hand side acting in the variable $x$.

\begin{theorem}
\label{t.cCf}
Let $(E,H)$ be an electromagnetic wave in $\mathcal{X}$ continuous up to
   $\overline{\mathcal{S}}$.
Then
$$
   \binom {E (x)}{H (x)}
 = \lim_{N \to \infty}
   \int_{\mathcal{S}}
   R_N (x,\cdot)
   \binom {\imath n (H)}
          {- \nu \wedge t (E)}
   ds
$$
for all $x \in \mathcal{X} \setminus \partial \mathcal{X}$.
\end{theorem}

In \cite{Tark95}, several approaches to the Cauchy problem for solutions of
linear elliptic equations with data on a part of the boundary are elaborated.
They give not only explicit formulas of Carleman type for solutions but also
conditions on the Cauchy data which are necessary and sufficient for the
Cauchy problem to be solvable.
The results demonstrate rather strikingly that the Cauchy problem for elliptic
equations is overdetermined.
It is solvable for a thin set of Cauchy data, and so any error in the Cauchy
data leads to unsolvability.
For this reason the variational approach to the Cauchy problem for elliptic
equations is of practical interest.
It works also in the case of nonlinear elliptic equations, cf.
   \cite{LyTark09}.
More generally, the problem consists in introducing reasonable classes of
approximate solutions.
In this sense any Carleman formula like the formula of Theorem \ref{t.cCf}
produces approximate solutions to the Cauchy problem with data on a boundary
piece for Maxwell's equations.

\end{document}